\documentclass[12pt]{article}
\usepackage{amsfonts}
\usepackage{amsmath}
\usepackage{amssymb,amscd,amsthm,array}

\let\ssection=\section
\renewcommand{\section}{\setcounter{equation}{0}\ssection}
 
\newcommand{\gl}{{\mathrm{gl}}}
\newcommand{\bbR}{\mathbb{R}}
\newcommand{\bbC}{\mathbb{C}}
\newtheorem{df}{Definition}[section]
\newtheorem{thm}{Theorem}[section]
\newtheorem{lem}{Lemma}[section]

\title{Deformations of the Heisenberg algebra\\ inside $\gl(3,\mathbb{K})$}
\author {Ya\"el Fr\'egier }  

\newcommand{\fh}{\mathfrak{h}}
\newcommand{\fg}{\mathfrak{g}}

\begin{document}
\date{}
\maketitle

\begin{abstract}
We study non-trivial deformations of the natural imbedding of the Lie
algebra $\fh_1$ of lower triangular matrices (the Heisenberg Lie
algebra) into $\gl(3,\mathbb{K})$, where $\mathbb{K}=\bbR$ or $\bbC$. 
Our first result is the calculation of the first cohomology space $H^1(\fh_1;\gl(3,\mathbb{K}))$. We prove that
there are no obstructions for integrability of infinitesimal deformation and, furthermore, give an explicit
formula for the most general deformation.
\end{abstract}

\section{Introduction} 

The  Heisenberg algebra $\fh_1$ is the Lie algebra of three by three lower triangular matrices.
Our aim is to deform the standard imbedding 
\begin{equation}
\label{Stand}
\rho:\fh_1\hookrightarrow\gl(3,\mathbb{K}),
\end{equation}
where
$\mathbb{K}=\mathbb{R}$ or $\mathbb{C}$. 
The theory of deformations of Lie stuctures (algebras and morphisms) is now a classical
subject (see e.g. \cite{GER}, \cite{N-R}, \cite{RIC}). However, a new concept of miniversal
deformation of Lie algebras has been introduced in \cite{FF1}. It is of course inspired by the notion of
universal unfolding in singularity theory. While Fialowski and Fuchs developed this notion for Lie algebras, the
case of Lie algebra homomorphisms has been recently considered (see~\cite{AA1} and also \cite{OR1,OR2}).

In this paper we will completely describe the miniversal deformation of the imbedding (\ref{Stand}).

\section{Deformations of homomorphisms}

Let $\rho:\fh\to\fg$ be a homomorphism of Lie algebras. A deformation of $\rho$ is an expression of the type
\begin{equation} 
\label{sum}
\rho(t)=
\rho_{0} + 
\sum_{m=1}^\infty \rho_{m}(t)
\end{equation}
where $t=(t_1,\dots ,t_r)$ are the parameters of the deformation and
each term $\rho_{m}(t)$ is a linear map from $\fh$  to $\fg$ homogeneous in $t$ of degree $m$. Note that
$t_1,\dots ,t_r$ can be considered as real (or complex) parameters, or as generators of a commutative
associative algebra (cf. \cite{FF1}), upon the context.

The deformation must be a Lie homomorphism for every value of the parameter $t$, i.e. satisfy:
\begin{equation} 
\label{morp}
\rho(t)([X,Y])=[\rho(t)(X),\rho(t)(Y)] 
\end{equation}
for every $X,Y \in \fh$.

If the sum (\ref{sum}) is finite, the deformation is said to be polynomial. In our case, one can restrict to
polynomial deformations since our construction will involve only a finite number of non vanishing terms (cf. main
results below).

\subsection{Equivalent deformations and the first cohomology}

The standard Chevalley-Eilenberg differential is given, in the case of a linear map $m$ from $\fh$  to
$\fg$, by the following formula:
\begin{equation} 
\label{DiffOper}
\delta ^1 m(X,Y)=m([X,Y])-[\rho(X),m(Y)]+[\rho (Y),m(X)] .
\end{equation}
Let us expand formula (\ref{morp}) as a series in $t$, the first order term is of the form
\begin{equation} 
\label{inf}
\rho_1(t)=\sum_{i=1}^rt_i\,\rho_1^i.
\end{equation}
>From (\ref{morp}) one obtains $\delta\rho_1^i =0$, that is, each map $\rho_1^i$ is a one-cocycle.

Two deformations $\rho (t)$ and $\rho '(t)$
are equivalent if there exists an inner automorphism
$I(t):\fg\otimes\mathbb{K}[t_1,\dots ,t_r]\longrightarrow \fg\otimes\mathbb{K}[t_1,\dots ,t_r]$ of the form
$$ 
I(t)=
\exp\left(\sum_{1\leq i\leq r}t_i\,ad A_i+\sum_{1\leq i,j\leq r}t_it_j\,ad A_{ij}+\dots\right),
$$ where
$A_{i},A_{ij},\dots$ are some elements of $\fg$, such that the relation
$
I(t)\circ\rho (t)=\rho '(t)
$
is satisfied.

Moreover one can check that the first order terms  $\rho_1^i$ 
and ${\rho '}_1^i$ differ by a coboundary i.e. ${\rho '}_1^i=\rho_1^i+\delta A_i$.
It follows that infinitesimal deformation of the homomorphism $\rho$ are classified by the first cohomology
space $H^1(\fh;\fg)$.

\subsection{Cup-product and Maurer-Cartan equation}

The standard cup-product (or Nijenhuis-Richardson product) of linear maps 
$a,b:\fh\to\fg$ is the linear map $[\![ a ,b ]\!] : \fh\otimes\fh\to\fg$ defined by
$$
[\![ a ,b ]\!](X,Y)=[a(X),b(Y)] - [a(Y),b(X)].
$$

Put $\phi (t)=\rho (t) -\rho_{0}$, the morphism equation (\ref{morp}) reads
\begin{equation}
\label{MC}
\delta \phi (t)-\frac{1}{2} [\![ \phi (t),\phi (t) ]\!]=0
\end{equation}
(see \cite{N-R,RIC}). This equation is called the Maurer-Cartan
equation (or the deformation equation).

Developing the Maurer-Cartan equation (\ref{MC}), one gets
\begin{equation} 
\label{MCOrdm}
\delta\rho_{m}(t)
=\frac{1}{2} \sum _{i+j=m} [\![ \rho _{i}(t),\rho _{j}(t) ]\!]
\end{equation}
for each $m$. 
The right hand side of this equation is always a 2-cocycle for any $m$ (cf., e.g., \cite{Fu}).
The equation admits a solution if and only if it is a coboundary. The cohomology class of the 2-cocycle in the
right hand side of~(\ref{MCOrdm}) is an obstruction for prolongation of the deformation to the order $m$.

\subsection{Construction of the miniversal deformation}\label{alg}

We are interested in deformations up to equivalence, hence, we will set $r=\dim H^1(\fh;\fg)$ and choose the
basis $[c_1],\ldots,[c_r]$ of $H^1(\fh;\fg)$ where $c_1,\ldots,c_r$ are non-trivial 1-cocycles on $\fh$ with
coefficients in $\fg$. We then put 
$$\rho_1^i=c_i.$$

The construction of the miniversal deformation goes as follows. Assume, by induction, that we constructed the
deformation (\ref{sum}) to the order $m-1$. To construct the $m$-th order term, one has to solve the equation
(\ref{MCOrdm}).
The right hand side of~(\ref{MCOrdm}) is an element of $Z^2(\fh;\fg)\otimes\mathbb{K}_m[t_1,\ldots,t_r]$; if
this is a coboundary then there exists a solution of (\ref{MCOrdm}). 

The solution of~(\ref{MCOrdm}) can be chosen arbitrarily up to the equivalence and reparametrization. Indeed, if
$\rho_{m}(t)$ and $\rho'_{m}(t)$ are two solutions, then their difference is a 1-cocycle

\section{The main results}

We formulate here the main results of this paper, all proofs will be given in Section \ref{ProofSec}.

\subsection{The first group of cohomology.}

One determines $H^1(\fh_1, \gl(3,\mathbb{K}))$ in order to know the
dimension of the parameter space and the infinitesimal generators of the deformation.

\begin{thm}\label{thm1}
Dim $H^1(\fh_1, \gl(3,\mathbb{K}))$=4.
\end{thm}

We also give a basis of this cohomology space. 
Let $e_{ij} \in \gl(3,\mathbb{K})$ be the standard basis of $\gl(3,\mathbb{K})$, namely $(e_{ij})_{kl} =
 \delta_{ik} \delta_{jl}$. Denote $\mathfrak{B}$ the natural basis of $\fh_1 $: 
$$
X  = e_{21},
\qquad 
Y  = e_ {31} ,
\qquad 
Z  =
 e_{32}
$$
and let $X^*,Y^*,Z^*$ be the dual basis in $\fh_1^*$.
The non-trivial cohomology classes are generated by the following
1-cocycles.
\begin{equation}
\label{OneCoc}
\begin{array}{l}
\rho^1  =  X^*\otimes e_{32} \\[4pt]
\rho^2  =  Z^*\otimes e_{21} \\[4pt]
\rho^3  =  Z^*\otimes(2 e_{22}+ e_{11})+Y^*\otimes e_{21}\\[4pt]
\rho^4  =  X^*\otimes( e_{11}+ e_{22}+ e_{33}).
\end{array}
\end{equation}


\subsection{An expression for the miniversal deformation.}

We will apply the algorithm described in Section \ref{alg} to the 1-cocycles (\ref{OneCoc}). The result is an explicit
formula for the miniversal deformation.

\begin{thm}
\label{thm2}
Up to equivalence and reparametrisation, the miniversal deformation of $\rho$ is given by the formula
\begin{eqnarray*}
\widetilde{\rho}(t_1,t_2,t_3,t_4) & = & \rho+t_1\,\rho^1+t_2\,\rho^2+t_3\,\rho^3+t_4\,\rho^4\\[6pt]
                                  && +  t_1t_2\,\rho^{12}+t_1t_3\,\rho^{13}+t_3^2\,\rho^{33}+t_1t_3^2\,\rho^{133}
\end{eqnarray*}
where
$$
\begin{array}{l}
\rho^{12}=X^*\otimes e_{21}\\[4pt]
\rho^{13}=X^*\otimes(2 e_{22}+ e_{11})\\[4pt]
\rho^{33}=-2Z^*\otimes e_{23}\\[4pt]
\rho^{133}=-2X^*\otimes e_{23}.
\end{array}.
$$

Applied on a generical element this gives:
 $$ \widetilde{\rho} (t_1,t_2,t_3,t_4)
\left(
 \begin{array}{ccc}
0     &0        &0  \\

a     &0        &0 \\

c     &b        &0 \\

 \end{array}
\right) = \left( 
 \begin{array}{ccc}
bt_3+at_1t_3+at_4     &0                    &0  \\

a+bt_2+at_1t_2+ct_3   &2bt_3+2at_1t_3+at_4  &-bt_3^2-at_1t_3^2\\

c                     &b+t_1                &at_4\\

 \end{array}
)\right).$$

\end{thm}

\noindent
Note, in comparison with \cite{AA1}, that no obstruction to the integrability appears, and
hence the parameter space is the free commutative algebra $\mathbb{K}[t_1,t_2,t_3,t_4]$. In other words, there
are no relations on the parameters of the deformation.

\section{Proofs of the main results}\label{ProofSec}

\subsection{Proof of theorem \ref{thm1}}

The proof consists in two steps. We will first calculate the space of 1-cocycles and then determine its subspace of
coboundaries.

\subsubsection{Computation of $Z^1(\fh_1, \gl(3,\mathbb{K}))$}
   
Let us calculate explicitly the expression of the differential $\delta ^1$ defined by the formula (\ref{DiffOper}). 
We will then determine its kernel.
    
\begin{lem}
\label{BasLem}
A basis of $Z^1$ is given by the vectors:
\begin{eqnarray*}
 e_1 & = & X^*\otimes  e_{32},\\
 e_2 & = & Z^*\otimes  e_{21},\\
e_3 & = & X^*\otimes  e_{33}+X^*\otimes  e_{22}+X^*\otimes  e_{11},\\
 e_4 & = & Z^*\otimes  e_{22}+\frac{1}{2}Z^*\otimes  e_{11}+\frac{1}{2}Y^*\otimes  e_{21},\\
 e_5 & = & Z^*\otimes  e_{31},\\
e_6 & = & Y^*\otimes  e_{31}+X^*\otimes  e_{21},\\
e_7 & = &  Y^*\otimes  e_{32}+X^*\otimes  e_{22}-X^*\otimes  e_{11},\\
e_8 & = & X^*\otimes  e_{31},\\
e_9 & = & Z^*\otimes  e_{32}+Y^*\otimes  e_{31}, \\
e_{10} & = & X^*\otimes  e_{23}+Y^*\otimes ( e_{33}- e_{11})-Z^*\otimes  e_{12}, \\
e_{11} & = & Z^*\otimes  e_{33}+\frac{1}{2}Z^*\otimes  e_{11}-\frac{1}{2}Y^*\otimes  e_{21}.\\ 
\end{eqnarray*}
\end{lem}

  Let $\Phi$ be a one cochain on $\fh_1$ with coefficients in $\gl(3,\mathbb{K})$. 
Using the basis $X,Y,Z$ of $\fh_1$, let us write this cochain in the form:
$$
(\underbrace{\Phi _{X}^{11},\dots ,\Phi _{X}^{33}}_{\Phi (X)},\,\underbrace{\Phi 
_{Y}^{11},\dots ,\Phi _{Y}^{33}}_{\Phi (Y)},\,\underbrace{\Phi _{Z}^{11},\dots 
,\Phi _{Z}^{33}}_{\Phi (Z)})
$$ 
and express the definition of the coboundary operator $\delta ^1$. One has
\begin {eqnarray*}
  \delta ^1\Phi (X,Y)   & =  & \Phi 
   ([ e_{21}, e_{31}])  -  [ e_{21},\Phi (Y)] + [ e_{31},\Phi (X)] \\
                             & = &  - [ e_{21},\sum _{i,j=1}^3\Phi _Y^{ij}  e_{ij}] + 
      [ e_{31},\sum _{i,j=1}^3\Phi _X^{ij}  e_{ij}]    \\
                             & = &    - \sum _{i,j=1}^3\Phi _Y^{ij}\,[ e_{21},  e_{ij}] + 
      \sum _{i,j=1}^3\Phi _X^{ij}\,[ e_{31}, e_ {ij}]
\end {eqnarray*}
which gives
$$
\delta ^1\Phi (X,Y)  =  - \sum _{j=1}^3\Phi _Y^{1j} e_{2j}  + \sum _{i=1}^3\Phi 
       _Y^{i2} e_{i1}    +\sum _{j=1}^3\Phi _X^{1j} e_{3j}  - \sum _{i=1}^3\Phi 
       _X^{i3} e_{i1}.
$$
Similarly,
 $$
 \delta ^1\Phi (X,Z)  =  -\sum _{i,j=1}^3\Phi _{Y}^{ij}  e_{ij}   - \sum _{j=1}^3\Phi _{Z}^{1j} e_{2j}  
       + \sum _{i=1}^3\Phi 
       _{Z}^{i2} e_{i1}       +\sum _{j=1}^3\Phi _{X}^{2j} e_{3j}  - \sum _{i=1}^3\Phi 
       _{X}^{i3} e_{i2}  
$$
and
$$
 \delta ^1\Phi (Y,Z)   =     - \sum _{j=1}^3\Phi _{Z}^{1j} e_{3j}  + \sum _{i=1}^3\Phi 
       _{Z}^{i3} e_{i1}       +\sum _{j=1}^3\Phi _{Y}^{2j} e_{3j}  - \sum _{i=1}^3\Phi 
       _{Y}^{i3} e_{i2}. 
$$

Let us comute $\delta ^1$ as a matrix. Since $\delta^1\Phi$ is a 2 cochain, it can be written
in the basis
$$
X^\star\wedge Y^\star\otimes  e_{ij},
\quad
X^\star \wedge Z^\star \otimes
e_{ij},
\quad
Y^\star \wedge  Z^\star \otimes e_{ij}.
$$
More precilely,
$$
\delta^1\Phi=
\sum _{i,j=1}^3
\left(
(\delta^1\Phi) _{X,Y}^{ij}\, 
X^\star\wedge Y^\star+
(\delta^1\Phi) _{X,Z}^{ij}\,X^\star \wedge Z^\star
+(\delta^1\Phi) _{Y,Z}^{ij}\,Y^\star \wedge  Z^\star\right)\otimes e_{ij}
$$
Applying $\delta^1\Phi$ to $(X,Y)$ one gets 
$$
\delta ^{1}\Phi (X,Y) =\sum _{i,j=1}^3(\delta^1\Phi) _{X,Y}^{ij}\,
e_{ij}.
$$
One can then identify the first nine coefficients $(\delta^1\Phi)_{X,Y}^{ij} , 1\leq i,j \leq 3$ which
correspond to the matrix of $ \delta^1 $ first nine rows.
Applying the same procedure to $(\delta^1\Phi)_{X,Z}^{ij}$ and $(\delta^1\Phi)_{Y,Z}^{ij}$, one finally gets a
$(27\times 27)$-matrix, see Appendix \ref{a1}.
In order to determine  the kernel of $\delta^1$, one has to find a maximal free subfamily among the column
vectors of the matrix of $\delta^1$. Dependance relations among remaining vectors will then give the kernel. Details of
these computations can also be found in Appendix \ref{a1}. This completes the proof of Lemma \ref{BasLem}.

\subsubsection{Computation of $B^1(\fh_1,\gl(3,\mathbb{K}))$}

The space of coboundaries $B^1(\fh_1,\gl(3,\mathbb{K}))$ is the image of 
$\gl(3,\mathbb{K})$ by the operator  
$ \delta ^0:\bigwedge^0\longrightarrow 
\bigwedge^1 $ defined by 
$$
\delta ^0 (A )(a)=[A ,a ]
$$
where $a\in\fh_1$ and $A \in  \gl(3,\mathbb{K}) $.

Proceding as above, one has for $A = (A ^{ij})$ with $i,j 
={1,2,3}$ 
$$
\begin{array}{rcl}
   \delta ^{0}A (X)   & = &
\displaystyle  - [ e_{21},\sum _{i,j=1}^3A ^{ij}  e_{ij}] \\
      & = & 
\displaystyle - \sum _{j=1}^3A ^{1j} e_{2j}  + \sum _{i=1}^3A ^{i2} e_{i1}  \\
   \delta ^{0}A (Y)   & = &  
\displaystyle- \sum _{j=1}^3A ^{1j} e_{3j}  + \sum _{i=1}^3A ^{i3} e_{i1}  \\
   \delta ^{0}A (Z)   & = &   
\displaystyle- \sum _{j=1}^3A ^{2j} e_{3j}  + \sum _{i=1}^3A 
      ^{i3} e_{i2}.
\end{array}
$$
The matrix of this operator is given in Appendix \ref{a2}. A basis of the image is as follows.
\begin {eqnarray*}
\delta_{11} & = & -X^\star\otimes{21}-Y^\star\otimes{31},\\
\delta_{12} & = & X^\star\otimes({11}-{22})-Y^\star\otimes{32},\\
\delta_{13} & = & -X^\star\otimes{23}+Y^\star\otimes({11}-{33})+Z^\star\otimes{12},\\
\delta_{21} & = & -Z^\star\otimes{31},\\
\delta_{22} & = & X^\star\otimes{21}-Z^\star\otimes{32},\\
\delta_{23} & = & Y^\star\otimes{21}+Z^\star\otimes({22}-{33}),\\
\delta_{32} & = & X^\star\otimes{31}.\\
\end {eqnarray*}

\subsubsection{Computation of $H^1(\fh_1,\gl(3,\mathbb{K}))$}

The dimension of $Z^1(\fh_1,\gl(3,\mathbb{K}))$ is 11, the one of  $B^1(\fh_1,\gl(3,\mathbb{K}))$ is  7. Hence the
quotient space  $H^1(\fh_1,\gl(3,\mathbb{K}))$ has dimension $11-7=4$.
One can check, see Appendices \ref{a2}-\ref{a3}, that the first four elements $e_1,e_2,e_3,e_4$ are independent modulo
$B^1(\fh_1,\gl(3,\mathbb{K}))$. Their classes form a basis of $H^1$. Theorem \ref{thm1} is proved.

\subsection{Proof of Theorem \ref{thm2}}

We will show that every infinitesimal deformation is integrable. In other words, there are no obtsructions to
integrability.

\subsubsection{Intergability at order 2}

One needs to evaluate the cup-products $[\![\rho^i,\rho^j ]\!], 1\leq
i, j \leq 4$ of the cocycles~(\ref{OneCoc}). 
It turns out that the only non vanishing terms are three coboundaries:
$$
\begin{array}{rcl}[\![\rho^1 ,\rho^2 ]\!]  
&= &  
X^*\wedge Z^*\otimes e_{31}=\delta^1(X^\star\otimes e_{21})\\[4pt]
[\![\rho^1 ,\rho^3 ]\!] & = & 2 X^*\wedge Z^* \otimes e_{32} + X^*\wedge Y^*\otimes e_{31}=\delta^1(X^\star\otimes(2
e_{22}+ e_{11}))\\[4pt]
[\![\rho^3 ,\rho^3 ]\!]  & = & 2 Z^* \wedge Y^*\otimes e_{21}=\delta^1(-2Z^\star\otimes e_{23})
\end{array}
$$
Let us choose
$$
\begin{array}{l}
\rho^{12}  = X^\star\otimes e_{21}\\[4pt]
\rho^{13}  =X^\star\otimes(2 e_{22}+ e_{11})\\[4pt]
\rho^{33}  =-Z^\star\otimes e_{23} 
\end{array}
$$ 
and $\rho^{ij}  = 0$ otherwise. One extends the deformation to the order 2. The equation (\ref{MCOrdm})
is satisfied to the second order.

\subsubsection{Intergability at order 3}

At the third order, the only non vanishing products are:
$$
\begin{array}{rcl}
[\![\rho^{12} ,\rho^3 ]\!] 
&=&
\frac{1}{2}\, X^*\wedge   Z^*\otimes ( - e_{21})\\[4pt]
{[\![}\rho^{13} ,\rho^2 {]\!] }
&=&
\frac{1}{2}\,X^*\wedge   Z^*\otimes (  e_{21})\\[4pt]
{[\![}\rho^{13} ,\rho^3 {]\!]}
&=&
\frac{1}{2}\,X^*\wedge   Y^*\otimes (  e_{21})\\[4pt]
{[\![}\rho^{33} ,\rho^1 {]\!]}
&=& 
- Z^*\wedge   X^*\otimes  ( e_{22}- e_{33})
\end{array}
$$
Hence the only r-h-s expressions that might appear in (\ref{MCOrdm}) are
$$
\begin{array}{rcl}
[\![\rho^{33} ,\rho^1 ]\!] +[\![\rho^3 ,\rho^{31} ]\!]
&=&
X^\star\wedge Y^\star\otimes  e_{21}+X^\star\wedge Z^\star\otimes ( e_{22}- e_{33})\\[4pt]
&=&\delta (-X^\star\otimes  e_{23})
\end{array}
$$
and
$$
\begin{array}{rcl}
[\![\rho^{12} ,\rho^3 ]\!] +\underbrace{[\![\rho^1 ,\rho^{23} ]\!]}_{=0}+[\![\rho^2 ,\rho^{13} ]\!]
&=&
X^\star\wedge Z^\star\otimes (- e_{21})+X^\star\wedge Z^\star\otimes  e_{21}\\[4pt]
&=&0.
\end{array}
$$

This shows that that the deformation can be extended to the third order (there are no obstructions).
More precisely, we set  $\rho^{331}=-X^\star\otimes  e_{23}$ and 
$\rho^{ijk}=0$ if  $(i,j,k)\neq(3,3,1)$.

\subsubsection{Intergability at order 4}
\label{4ord}

  At the fourth order, the only non vanishing cup-products appearing in (\ref{MCOrdm}) are:

$$
\begin{array}{rcl}
{[\![}\rho^{133} ,\rho^3 {]\!]}
&=&X^*\wedge Z^*\otimes e_{23}\\[4pt]
{[\![}\rho^{13} ,\rho^{33} {]\!]}
&=&
-X^*\wedge Z^*\otimes e_{23}.
\end{array}
$$
Hence the only equation for which the r-h-s is not identicaly zero a-priori is~$\rho^{1333}$.
However, one has
$$
[\![\rho^{133} ,\rho^3 ]\!]+[\![\rho^{13} ,\rho^{33} ]\!]
=  -X^*\wedge Z^*\otimes e_{23}+X^*\wedge Z^*\otimes e_{23} =0,
$$
therefore no obstruction appears. Hence the deformation can be extended to the order 4, setting
$\rho^{ijkl}=0$ for all  $i,j,k,l$.

\subsubsection{Intergability at order 5}

At the fifth order, all the cup-products are zero. Indeed,  a straightforward computation shows that all the
terms of the type $[\![\rho^{ijk} ,\rho^{lm} ]\!]$ vanish.
Hence all the terms appearing in the general formula (\ref{MCOrdm}) vanish. There are no obstructions to
integrabilitiy, and one can extend the deformation to the fifth order setting $\rho^{ijklm}=0$ for all  $i,j,k,l,m$.

\subsubsection{Intergability at any order}

 At order 6, all the cup-products constituted by one term of order
 one and an other term of order five vanish since we have just set the fifth order terms to be zero.
It is just the same for the terms constituted by an element of order four and an other one of order two since we've
set the elements of order four to be zero (cf \ref{4ord}).
 At last the only cup-product constituted with two non vanishing terms of order three is 
$[\![\rho^{331} ,\rho^{331}]\!]$ wich turns out to be zero.

The analogous arguments are valid at every order. Theorem \ref{thm2} is proved.

\section {Appendix}
\subsection {Computations appearing in  proof of Theorem \ref{thm1}}

\subsubsection {Computing $\delta^1$}

The first fourteen elements of this familly are independent (indeed, each these element has an
$\underline{underlined}$ non vanishing component which vanishes for the other elements).
\label{a1}
Setting    
$$
\delta ^X_{ij}= \delta^1(X^\star \otimes e_{ij}),
\delta ^Y_{ij}= \delta^1(Y^\star \otimes e_{ij}),\delta
^Z_{ij}= \delta^1(Z^\star \otimes e_{ij}),
$$
we have
\begin{eqnarray*} 
\delta ^X_{11} & = & \underline{X^\star\wedge Y^\star\otimes  e_{31}},\\
\delta ^X_{12} & = & \underline{X^\star\wedge Y^\star\otimes  e_{32}},\\
\delta ^X_{13} & = & X^\star\wedge Y^\star\otimes (\underline{ e_{33}}-  e_{11})-X^\star\wedge Z^\star\otimes  e_{12},\\
\delta ^X_{21} & = & \underline{X^\star\wedge Z^\star\otimes  e_{31}},\\
\delta ^X_{22} & = & \underline{X^\star\wedge Z^\star\otimes  e_{32}},\\
\delta ^X_{23} & = & -X^\star\wedge Y^\star\otimes  e_{21}+X^\star\wedge Z^\star\otimes ( \underline{ e_{33}}- e_{22}),\\
 \delta ^Y_{11} & = & -X^\star\wedge Y^\star\otimes  e_{21}-\underline{X^\star\wedge Z^\star\otimes  e_{11}},\\
\delta ^Y_{12} & = & X^\star\wedge Y^\star\otimes ( e_{11}-\underline{ e_{22}})-X^\star\wedge Z^\star\otimes  e_{12},\\
\delta ^Y_{13} & = & -\underline{X^\star\wedge Y^\star\otimes  e_{23}}-X^\star\wedge Z^\star\otimes  e_{13}-Y^\star\wedge Z^\star\otimes  e_{12},\\
\delta ^Y_{21} & = & -\underline{X^\star\wedge Z^\star\otimes  e_{21}}+Y^\star\wedge Z^\star\otimes  e_{31},\\
\delta ^Y_{22} & = & X^\star\wedge Y^\star\otimes  e_{21}-X^\star\wedge Z^\star\otimes  e_{22}+\underline{Y^\star\wedge Z^\star\otimes  e_{32}},\\
\delta ^Y_{23} & = & -X^\star\wedge Z^\star\otimes  e_{23}+Y^\star\wedge Z^\star\otimes ( e_{33}-\underline{ e_{22}}),\\
\delta ^Z_{13} & = & -X^\star\wedge Z^\star\otimes  e_{23}-Y^\star\wedge Z^\star\otimes ( e_{33}-\underline{ e_{11}}),\\
\delta ^Z_{23} & = & \underline{Y^\star\wedge Z^\star\otimes
  e_{21}},
\end{eqnarray*}
One can check that the two following terms are indepent from the preceding ones.
\begin{eqnarray*}
\delta ^Y_{33} & = & -X^\star\wedge Z^\star\otimes  e_{33}-Y^\star\wedge Z^\star\otimes  e_{32}.\\
\delta ^Z_{11} & = & -X^\star\wedge Z^\star\otimes
e_{21}-Y^\star\wedge Z^\star\otimes  e_{31},
\end{eqnarray*}
The remaining elements 
\begin{eqnarray*}
\delta ^X_{31} & = & 0,\\
\delta ^X_{32} & = & 0,\\
\delta ^X_{33} & = & -X^\star\wedge Y^\star\otimes  e_{31}-X^\star\wedge Z^\star\otimes  e_{32}.\\
\delta ^Y_{31} & = & -X^\star\wedge Z^\star\otimes  e_{31},\\
\delta ^Y_{32} & = & X^\star\wedge Y^\star\otimes  e_{31}-X^\star\wedge Z^\star\otimes  e_{32},\\
\delta ^Z_{12} & = & -X^\star\wedge Z^\star\otimes ( e_{22}- e_{11})-Y^\star\wedge Z^\star\otimes  e_{32},\\
\delta ^Z_{21} & = & 0,\\
\delta ^Z_{22} & = & X^\star\wedge Z^\star\otimes  e_{21},\\
\delta ^Z_{31} & = & 0,\\
\delta ^Z_{32} & = & X^\star\wedge Z^\star\otimes  e_{31},\\
\delta ^Z_{33} & = & Y^\star\wedge Z^\star\otimes  e_{31}
\end{eqnarray*}
are linear combinations of the previous ones. More precisely,
$\delta ^X_{31}$,$\delta ^X_{32}$,$\delta ^Z_{21}$ and $\delta ^Z_{31}$ vanish and
\begin{eqnarray*}
\delta ^X_{33} & = & -\delta ^X_{22}-\delta ^X_{11}\\
 \delta ^Y_{31} & = & -\delta ^X_{21}\\
 \delta ^Y_{32} & = & -\delta ^X_{22}+\delta ^X_{11}\\
\delta ^Z_{22} & = & -\frac{1}{2}\delta ^Z_{11}-\frac{1}{2}\delta ^Y_
{21}\\
  \delta ^Z_{32} & = & -\delta ^Y_{31} \\
  \delta ^Z_{12} & = & -\delta ^Y_{11}+\delta ^X_{23}+ \delta ^Y_{33} \\
\delta ^Z_{33} & = & -\frac{1}{2}\delta ^Z_{11}+\frac{1}{2}\delta
 ^Y_{21}.\\ 
 \end{eqnarray*}

These relations are important since they give us the basis of $Z^1$:
\begin{eqnarray*}
 e_1 & = & \underline{X^*\otimes  e_{32}},\\
 e_2 & = & \underline{Z^*\otimes  e_{21}},\\
 e_3 & = & \underline{X^*\otimes  e_{33}}+X^*\otimes  e_{22}+X^*\otimes  e_{11},\\
 e_4 & = & Z^*\otimes  e_{22}+\frac{1}{2}\underline{Z^*\otimes  e_{11}}+\frac{1}{2}Y^*\otimes  e_{21},\\
e_5 & = & Z^*\otimes  e_{31},\\
e_6 & = & Y^*\otimes  e_{31}+X^*\otimes  e_{21},\\
e_7 & =  & Y^*\otimes  e_{32}+X^*\otimes  e_{22}-X^*\otimes  e_{11},\\
e_8 & = & X^*\otimes  e_{31},\\
e_9 & = & Z^*\otimes  e_{32}+Y^*\otimes  e_{31}, \\
e_{10} & = & X^*\otimes  e_{23}+Y^*\otimes ( e_{33}- e_{11})-Z^*\otimes  e_{12}, \\
e_{11} & = & Z^*\otimes  e_{33}+\frac{1}{2}Z^*\otimes  e_{11}-\frac{1}{2}Y^*\otimes  e_{21}.\\ 
\end{eqnarray*}
These elements are in the kernel because of the preceding relations.

\subsubsection {Computing $\delta^0$}
\label{a2}

 The matrix of $\delta ^0$ is given by the following set of vectors:
\begin{eqnarray*}
\delta_{11} & = & -X^\star\otimes{21}-\underline{Y^\star\otimes{31}},\\
\delta_{12} & = & X^\star\otimes(\underline{{11}}-{22})-Y^\star\otimes{32},\\
\delta_{13} & = & -\underline{X^\star\otimes{23}}+Y^\star\otimes({11}-{33})+Z^\star\otimes{12},\\
\delta_{21} & = & -\underline{Z^\star\otimes{31}},\\
\delta_{22} & = & X^\star\otimes{21}-\underline{Z^\star\otimes{32}},\\
\delta_{23} & = & \underline{Y^\star\otimes{21}}+Z^\star\otimes({22}-{33}),\\
\delta_{32} & = & \underline{X^\star\otimes{31}},
\end{eqnarray*}
which are independent and
\begin{eqnarray*}
\delta_{31} & = & 0,\\
\delta_{33} & = & Y^\star\otimes{31}+Z^\star\otimes{32},
\end{eqnarray*}
which are linear combinations of the above ones. Indeed, $\delta_{31}=0$ and $\delta_{33}=-\delta_{11}+\delta_{22}$.

\subsubsection {Computing the basis of $H^1(\fh_1,\gl(3,\mathbb K))$}
\label{a3}
The first four elements $e_1, e_2, e_3, e_4$ of the basis of the space of cocycles $Z^1(\fh_1,\gl(3,\mathbb K))$ are
linearly independent modulo the coboundaries, since each of them has an $\underline{underlined}$ component which does
not appear in the elements of the basis of the space of coboundaries
$B^1(\fh_1,\gl(3,\mathbb K))$.

\vskip 1cm

\noindent
{\bf Acknowledgments}. 
The problem was stated to me by B.Konstant and V. Ovsienko; I am  grateful to V. Ovsienko for his help. I would
also like to thank C. Duval for enlightening discussions, C. Roger and F. Pellegrini and  for their
interest in this work.

\vskip 1cm


\noindent
Ya\"el FREGIER\\
{\small C.N.R.S., C.P.T.}\\
{\small  Luminy-Case 907}\\
{\small  F-13288 Marseille Cedex 9, France}\\
\&\\
{\small Institut
Girard Desargues, URA CNRS 746}
\\ {\small Universit\'e Claude Bernard - Lyon I}\\
{\small 43
bd. du 11 Novembre 1918}\\
{\small 69622 Villeurbanne Cedex, France}

\end{document}